\documentclass{amsart}

\usepackage{times}
\usepackage{amssymb}
\usepackage{amsmath}
\usepackage{amsthm}
\usepackage{latexsym}
\usepackage{epsf}
\usepackage{graphicx}
\usepackage{color}
\usepackage{setspace}
\usepackage{xypic}
\usepackage{xy}
\usepackage{textcomp}
\usepackage{tikz}
\usepackage{enumerate}

\newcommand{\qedbox}{\nobreak \ifvmode \relax \else
      \ifdim\lastskip<1.5em \hskip-\lastskip
      \hskip1.5em plus0em minus0.5em \fi \nobreak
      \vrule height0.75em width0.5em depth0.25em\fi}
\def\brc#1{\langle #1 \rangle}
\def\noin{\noindent}

\def\lcm#1{[#1]}
\def\dblbar{|\hspace{-.02in}{|}}
\def\strdiv{|_{_{<}}}
\newcommand*\circleme[1]{\tikz[baseline=(char.base)]{
            \node[shape=circle,draw,inner sep=1pt] (char) {#1};}}

\newenvironment{pf}[1][Proof:]{\begin{trivlist}
\item[\hskip \labelsep {\itshape #1}]}{\end{trivlist}}

\swapnumbers
\theoremstyle{plain}
\newtheorem{theorem}{Theorem}[section]

\newtheorem{proposition}[theorem]{Proposition}
\newtheorem{lemma}[theorem]{Lemma}
\newtheorem{corollary}[theorem]{Corollary}

\theoremstyle{definition}
\newtheorem{definition}[theorem]{Definition}

\newtheorem{example}[theorem]{Example}

\theoremstyle{remark}
\newtheorem{remark}[theorem]{Remark}

\begin{document}

\thispagestyle{empty}

\title[Behavior of Trivariate Monomial Ideals]{On the Behavior of Minimal Free Resolutions of Trivariate Generic Monomial Ideals}
\author{Jared L Painter}

\begin{abstract}
 We will explore some properties of minimal graded free resolutions of $R/I$, where $R$ is a trivariate polynomial ring over a field and $I$ is a monomial ideal.  Our focus will be to consider a specific form of the resolutions when $I$ is primary to the homogeneous maximal ideal.  We will identify certain characteristics of the last matrix of these resolutions, and observe differences in the resolutions for generic ideals in comparison to non-generic ideals.  Finally, we learn how to identify whether $I$ is generic by knowing the structure of the last matrix in the minimal free resolution of $R/I$.
\end{abstract}

\maketitle

\section{Introduction}

For the duration of this paper we will assume that $R=\Bbbk[x,y,z]$ is a trivariate polynomial ring over a field $\Bbbk$ and $I$ is a monomial ideal that is primary to the homogeneous maximal ideal $\mathfrak{m}$ of $R$.  The previous conditions imply that $R/I$ is in fact local and artinian, but we will assume that it is not Gorenstein.  We will be discussing minimal graded free resolutions of $R/I$. In this case it is known that the projective dimension of $R/I$ is 3, thus the minimal free resolution of $R/I$ has the form,
$$\mathbb{F}:= \ 0 \longrightarrow F_3 \stackrel{f_3}{\longrightarrow} F_2 \stackrel{f_2}{\longrightarrow} F_1 \stackrel{f_1}{\longrightarrow} F_0 \longrightarrow R/I \longrightarrow 0.$$
In \cite[Theorem 20.9]{DE} it is shown that if any complex of free $R$ modules is exact then rank$F_i = $ rank$f_i + $ rank$f_{i+1}$.  Using this formula along with the assumption that $I$ is minimally generated by $n$ elements and rank$F_3 = m$ we obtain a more precise construction of the free resolution,
$$\mathbb{F}:= \ 0 \longrightarrow R^m \stackrel{f_3}{\longrightarrow} R^{m+n-1} \stackrel{f_2}{\longrightarrow} R^n \stackrel{f_1}{\longrightarrow} R \longrightarrow R/I \longrightarrow 0 \textrm{ \ with } n>3, m>1.$$
In this paper we think of the maps $f_i$ as matrices with respect to the standard bases of the free modules $R^j$.  We are interested in when we get nonzero elements from $I$ as entries in $f_i$ and what is the maximum number of such entries we can get in $f_i$.  If $I$ is minimally generated by monomials $m_1, \ldots , m_n$ then the matrix $f_1$ is $[m_1 \ldots m_n]$.  Since $f_2$ is generated by the S-pairs between the minimal generators, defined in ~\ref{SecSyzDef}, we will only get nonzero elements from $I$ in $f_2$ when the gcd$(m_i,m_j) = 1$ and the S-pair between $m_i$ and $m_j$ is a minimal second syzygy.  From this we will find that if $I$ is minimally generated by $n$ monomials, then $2n-2$ is the maximum number of nonzero entries we can get in $f_2$ which are also elements of $I$.  Our primary focus will be on determining when we get elements from $I$ in $f_3$.  We will use a special form of the minimal free resolution, given in Definition~\ref{ResCons} to find the maximum number of rows of $f_3$ that contain only elements from $I$.  We are primarily interested in generic monomial ideals. When $I$ is generic the free resolution for $R/I$ will have a nice construction. More precisely, if $I$ is generic then each column of the matrix of $f_3$ in the minimal free resolution of $R/I$ will have exactly three nonzero \emph{pure power} entries.  In Theorem~\ref{MaxEntriesFromI} we will show that for generic monomial ideals the maximum number of nonzero elements from $I$ that we get in $f_3$ is $n-2$, when $I$ is minimally generated by $n$ elements.  We will conclude by considering some examples of resolutions for non-generic monomial ideals and contrast the differences of $f_3$ when $I$ is generic.  We will also consider examples where we relax the conditions that $I$ is a \emph{trivariate} monomial ideal and that $I$ is $\mathfrak{m}-$primary.  Many of our results will not hold without these conditions.

One of the motivations for this work relates to a question on the Bass numbers of $R/I$.  Specifically, is the first nonzero Bass number of $R/I$ always smaller than the second nonzero Bass number of $R/I$?  Since we are assuming that $R/I$ is local and Artinian, the first nonzero Bass number of $R/I$ will be the zeroth Bass number.  In this case it is also known that the Bass numbers of $R/I$ are equal to the Betti numbers of the canonical module of $R/I$, see \cite{BH} and \cite{DJL}.  This allows us to relate the question on when we get nonzero entries from $I$ in $f_3$ to the question on Bass numbers, by finding a free presentation of the canonical module of $R/I$ from the minimal free resolution of $R/I$.  We can then restate the question on Bass numbers to be, after permissable row operations, is the number of rows in $f_3$ which contain only entries from $I$ less than or equal to $n-2$?  We discuss this in more detail in Remark~\ref{BassRmk} and find that we get a positive answer to this question when $I$ is generic.

\section{Preliminaries}

  In this section we begin by introducing notation and describing our special form of the minimal free resolution of $R/I$.  This will be followed by a description of the minimal and non-minimal second syzygies of the resolution in Lemma~\ref{SyzLem} and Remark~\ref{LinCombTwo}.  We will conclude by discussing when we get nonzero elements of $I$ as entries in $f_2$ and prove in Proposition~\ref{MaxEntriesf2} that the maximum number of such entries will be $2n-2$.

  We will denote the least common multiple of monomials $m_1, \ldots, m_r$ by $m_{1\ldots r}$. In particular, $m_{ij}=[m_i,m_j]$ denotes the least common multiple of $m_i$ and $m_j$.  Throughout this paper the monomial $m_i$ will be represented by $x^{a_i}y^{b_i}z^{c_i}$.  We will also say that a monomial $m'$ \emph{strongly divides} a monomial $m$, denoted $m' \dblbar m$, if $m'$ divides $m/x_i$ for all variables $x_i$ dividing $m$.  If a monomial $m'$ strictly divides a monomial $m$ we will write $m'\strdiv m$.  The following lemma is a simple fact, which we will use frequently in this paper.

\begin{lemma}
\label{DividesLem}
Let $m_i,m_j,$ and $m_k$ be distinct minimal generators of $I$, then
\begin{enumerate}
\item $m_k|m_{ij}$ if and only if both $m_{ik}$ and $m_{jk}$ divide $m_{ij}$
\item $m_{ik}|m_{ij}$ if and only if $m_{jk}|m_{ij}$
\item if $m_k\dblbar m_{ij}$ then $m_{ik}\strdiv m_{ij}$ and $m_{jk} \strdiv m_{ij}$.
\end{enumerate}
\end{lemma}

\begin{pf}
To prove this we only need consider the exponents on one of the variables in these monomials.  It is easy to see we can extend our argument to the rest of the variables.

\noin (1): If $m_k|m_{ij}$ then $a_k \leq $ max$\{a_i,a_j\}$, which implies that max$\{a_i,a_k\}$ and max$\{a_j,a_k\}$ are both less than or equal to max$\{a_i,a_j\}$.  Thus both $m_{ik}$ and $m_{jk}$ divide $m_{ij}$.  The reverse direction is similar.

\noin (2): If $m_{ik}|m_{ij}$ then max$\{a_i,a_k\} \leq$ max$\{a_i,a_j\}$, which implies that max$\{a_j,a_k\} \leq$ max$\{a_i,a_j\}$, thus $m_{jk}|m_{ij}$.  The reverse direction is the same.

\noin (3): Since $m_i$ and $m_j$ are both minimal generators of $I$, then without loss of generality we may assume that $a_i>a_j$ and $b_j>b_i$. By definition, if $m_k\dblbar m_{ij}$, then $a_k < $ max$\{a_i,a_j\}=a_i$ and $b_k < $ max$\{b_i,b_j\}=b_j$.  This implies that max$\{a_j,a_k\} < $ max$\{a_i,a_j\}$ and max$\{b_i,b_k\} < $ max$\{b_i,b_j\}$.  Thus $m_{ik}\strdiv m_{ij}$ and $m_{jk} \strdiv m_{ij}$. \qed
\end{pf}

This lemma will help us describe the minimal second syzygies. Generally speaking the second syzygies are generated by the S-pairs of the minimal generators of $I$.  We will now define the generators for the second syzygies of $R/I$.
\begin{definition}
\label{SecSyzDef}
If $I$ is a monomial ideal with minimal generating set $\{m_1, \ldots, m_n\}$, then the second syzygies of $R/I$ are generated by the S-pairs between the minimal generators of $I$:
$$\displaystyle \sigma_{ij} = \frac{m_{ij}}{m_j}e_j - \frac{m_{ij}}{m_i}e_i, \textrm{ for } 1 \leq i < j \leq n.$$
\end{definition}
The set $\{\sigma_{ij}\}_{i<j}$ is rarely a minimal generating set for the second syzygies of $R/I$.  We denote the set of all second syzygies by $Z_2 = \sum_{i<j}R\sigma_{ij} \subseteq R^n$.  We will define a unique minimal generating set for $Z_2$ denoted $S_2$, which we will call the \textit{ordered minimal second syzygies}. This minimal generating set will be defined from the following orderings on the generators of $I$ and the $\sigma_{ij}$'s.  For the duration of this paper we will use the graded reverse lexicographic ordering (GRevLex) on $R$ with $x<y<z$, and define a standard dictionary order on the indices of $\sigma_{ij}$.  That is, $\sigma_{ij} < \sigma_{kl}$ if and only if either $i<k$, or $j<l$ when $i=k$.  It should be noted that when we write $m_{ij}$ it is not implied that $i<j$, since $m_{ij}=m_{ji}$.  However when we write $\sigma_{ij}$, it is always assumed that $i<j$.  We can now define a specific minimal generating set $S_2$ of $Z_2$, which we will use throughout this paper.

\begin{definition}
\label{ResCons}
If $I$ is a monomial ideal, then $S_2$ is the set of \textit{ordered minimal second syzygies} of $R/I$, such that $\sigma_{ij} \in S_2$ if and only if the following conditions are satisfied,
\begin{enumerate}\setlength{\itemsep}{-0pt}
\item $\sigma_{ij} \in Z_2-\mathfrak{m}Z_2$ and
\item $\sigma_{ij} \not= \sum_{k<l} a_{kl}\sigma_{kl}$, $a_{kl} \in R$, in which $\sigma_{ij} < \sigma_{kl}$ for all $k,l$ such that $a_{kl}$ is a unit.
\end{enumerate}
\end{definition}

By Nakayama's lemma the minimal generators of $Z_2$ will be contained in $Z_2-\mathfrak{m}Z_2$.  In some situations the problem arises where we may have a choice of which second syzygy we remove to construct a minimal generating set for $Z_2$.  If this occurs, then by Definition~\ref{ResCons} we will always remove the smallest second syzygy, with respect to GRevLex and the chosen dictionary ordering on $\{\sigma_{ij}\}_{i<j}$, of those from which we have a choice. In terms of the resolution of $R/I$, $S_2$ will consist of the $\sigma_{ij}$ which are the columns of the matrix of $f_2$.

In \cite{EM} and \cite{MS} we learn that the matrix for $f_3$ is \emph{completely determined} by the minimal second syzygies which are used in $f_2$.  This is because we can represent the minimal resolution of $R/I$ by some labeled planar graph.  The minimal second syzygies are represented by the edges of the planar graph and the minimal third syzygies are represented by the faces of the planar graph.  The labeling is given by the ordering chosen on the minimum generators of $I$ which are the vertices of the graph. The edges that are chosen will determine the faces of the graph.  These are the minimal cycles formed by the minimal second syzygies. These cycles are the faces of the planar graph which represents the minimal third syzygies.  For more on planar graphs we refer the reader to \cite{GR} and \cite{EM}.  Since we are choosing the maximal second syzygies based on the chosen ordering of $\{\sigma_{ij}\}_{i<j}$ we will refer to the resolution obtained from Definition~\ref{ResCons} as the \emph{maximal ordered resolution}. The resolutions we are considering here all have projective dimension 3, thus the column ordering of $f_3$ will not have an affect on any results of this paper.  To be consistent, we will choose the same ordering in $f_3$ that is being used in $f_2$.  The following lemma gives conditions on the minimal generators of $I$ so that an S-pair $\sigma_{ij} \not\in S_2$ and in turn gives us conditions for $\sigma_{ij}$ to be minimal.
\pagebreak
\begin{lemma}[Second Syzygy Lemma]
\label{SyzLem}
Let $I$ be a monomial ideal with minimal generating set $\{m_1, \ldots, m_n\}$.
\begin{enumerate}\setlength{\itemsep}{-0pt}
\item If $\sigma_{ij} \in \mathfrak{m}Z_2$ then there exists a minimal generator $m_k$ such that $m_{ik} \strdiv m_{ij}$.
\item If there exists a minimal generator $m_k$ such that $m_{ik}\strdiv m_{ij}$ and $m_{jk} \strdiv m_{ij}$ then $\sigma_{ij} \in \mathfrak{m}Z_2$.
\item If $\sigma_{ij} \not\in S_2$ then there exists a minimal generator $m_k$ such that $m_k\strdiv m_{ij}$.
\end{enumerate}
\end{lemma}

\begin{pf}
(1): If $\sigma_{ij} \in \mathfrak{m}Z_2$, then
$$\displaystyle \sigma_{ij} = \sum_{1\leq k < l \leq n} a_{kl}\sigma_{kl}, \textrm{ with } a_{kl} \in \mathfrak{m}.$$
It is enough to look at what we would need to get the $i^{\textrm{th}}$ row entry in $\sigma_{ij}$.  We have,
$$\displaystyle \frac{m_{ij}}{m_i} = \sum_{k=1}^{i-1} a_{ki}\frac{m_{ki}}{m_i} - \sum_{l=i+1}^n a_{il}\frac{m_{il}}{m_i} \ \Longrightarrow \ m_{ij} = \sum_{k=1}^{i-1} a_{ki}m_{ki} - \sum_{l=i+1}^n a_{il}m_{il}.$$
Since $m_{ij}$ is a monomial we know that at least one of the terms in the above sum must equal $m_{ij}$ up to multiplication by a unit.  That is, there exists a minimal generator $m_k$ such that $m_{ij} = ua_{ik}m_{ik}$ where $u$ is a unit in $R$.  Thus we have that $m_{ik} \strdiv m_{ij}$ since $a_{ik} \in \mathfrak{m}$.\\
(2): Since $m_{ik}\strdiv m_{ij}$ and $m_{jk} \strdiv m_{ij}$ we know that both $\displaystyle a = \frac{m_{ij}}{m_{ik}}, b = \frac{m_{ij}}{m_{jk}} \in \mathfrak{m}$.  Thus we have that,
$$a\sigma_{ik} - b\sigma_{jk} \ = \ \displaystyle a\left(\frac{m_{ik}}{m_k}e_k - \frac{m_{ik}}{m_i}e_i\right) - b\left(\frac{m_{jk}}{m_k}e_k - \frac{m_{jk}}{m_j}e_j\right)$$
$$\hspace{.1in} \displaystyle = \ \frac{m_{ij}}{m_k}e_k - \frac{m_{ij}}{m_i}e_i - \frac{m_{ij}}{m_k}e_k + \frac{m_{ij}}{m_j}e_j$$
$$ \displaystyle = \ \frac{m_{ij}}{m_j}e_j - \frac{m_{ij}}{m_i}e_i \ = \ \sigma_{ij}. \hspace{.52in}$$
This implies that $\sigma_{ij} \in \mathfrak{m}Z_2$.\\
(3): If the hypothesis of (1) from Definition~\ref{ResCons} is not satisfied, then $\sigma_{ij} \not\in S_2$. Thus from the proof of (1) we have that there is a minimal generator $m_k$ such that $m_k\strdiv m_{ij}$.  The only other possibility is that $$\sigma_{ij} = \sum_{k<l} a_{kl}\sigma_{kl}, a_{kl} \in R, \textrm{ in which } \sigma_{ij} < \sigma_{kl} \textrm{ for all } k,l \textrm{ such that } a_{kl} \textrm{ is a unit.}$$
Using a similar argument as in the proof of (1), we have that there exists a minimal generator $m_k$ such that $m_{ij} = ua_{ik}m_{ik}$, $a_{ik}$ is not necessarily in $\mathfrak{m}$.  Hence $m_{ik}|m_{ij}$ which implies that $m_k|m_{ij}$ by Lemma~\ref{DividesLem}.
Thus we must have that $m_k \strdiv m_{ij}$ since we are assuming that $m_k$ is a minimal generator of $I$.\qed
\end{pf}

\begin{remark}
\label{LinCombTwo}
There are two important facts that we obtain from the previous lemma. First, if $\sigma_{ij} \not\in S_2$ then we know that there must be a minimal generator $m_k$ of $I$ such that $a\sigma_{ik} + b\sigma_{jk} = \sigma_{ij}$. Thus every $\sigma_{ij} \not\in S_2$ can be obtained from a linear combination of exactly two second syzygies.  Secondly, the contrapositive of (3) says that if there is no minimal generator $m_k$ of $I$ such that $m_k\strdiv m_{ij}$ then $\sigma_{ij} \in S_2.$
\end{remark}

We can now address the question of when we get elements from $I$ as entries in $f_2$ and what the maximum number of such entries will be.  Since we are only dealing with trivariate monomial ideals, if two minimal generators $m_i$ and $m_j$ have gcd$(m_i,m_j) = 1$ then one of the generators must be a pure power and the other generator must only have nonzero degrees on the other two variables.  Using this in combination with Remark~\ref{LinCombTwo} we can classify the maximum number of nonzero entries of $f_2$ which are also in $I$.

\begin{proposition}
\label{MaxEntriesf2}
Let $I$ be a monomial ideal of $R$ minimally generated by $n$ elements.
\begin{enumerate}
\item If $m_i$ and $m_j$ are minimal generators of $I$ such that gcd$(m_i,m_j) =1$ and there is no other minimal generator $m_k$ of $I$ such that $m_k\strdiv m_{ij}$, then $$\sigma_{ij}\in S_2 \textrm{ and } \sigma_{ij} = m_ie_j - m_je_i.$$
\item The matrix of $f_2$ from Definition~\ref{ResCons} contains at most $2n-2$ nonzero entries from $I$.
\end{enumerate}
\end{proposition}
\begin{pf}
(1):  This follows immediately from Definition~\ref{SecSyzDef} and Remark~\ref{LinCombTwo}, which is a consequence of Lemma~\ref{SyzLem}.\\
(2):  Since we are assuming that $n>3$ we must have at least one generator with nonzero degrees on at least two variables.  For simplicity order the generators so that $m_1=x^a, m_2=y^b, m_3=z^c$.  Note we are not assuming that these are the only generators of $I$.  There are two cases that need to be considered.  First, if there are no minimal second syzygies between a pure power generator and a mixed double generator, then we can only get entries from $I$ in $f_2$ from the second syzygies $\sigma_{12}, \sigma_{13},$ and $\sigma_{23}.$  This gives a maximum of 3 columns in $f_2$ with entries from $I$, which means the number such columns is less than or equal to $n-1$ since $n\geq 4$, which satisfies our hypothesis using part (1) of the Proposition.

For the second case we will show that if we have a minimal second syzygy between a pure power generator and a mixed double generator involving the other two variables, then this is the only pure power generator that can appear more than once in $f_2$.  Without loss of generality assume that $x^a$ is the pure power generator we are interested in, and that $m_4=y^{\beta}z^{\gamma},$ $\beta,\gamma>0$ is the mixed double generator with $\sigma_{14}\in S_2$.  We will first show that the only second syzygies that could yield a $y^b$ or $z^c$ in $f_2$ must be $\sigma_{12}$ and $\sigma_{13}$ respectively.  This will imply that all of the second syzygies in $S_2$ that have nonzero entries from $I$ are of the form $\sigma_{1i}$ where $m_i = y^{b'}z^{c'}$.  Secondly we will show that the maximum number of these minimal second syzygies that we may have is $n-1$ when $I$ is minimally generated by $n$ elements. This implies that the maximum number of nonzero entries from $I$ that we can have in $f_2$ is $2n-2$ by part (1) of the Proposition.

It is implied that $\sigma_{23}\not\in S_2$ since $m_4$ is a minimal generator of $I$. Suppose that $\sigma_{25} \in S_2$ such that $m_5=x^{a'}z^{c'}, a',c'>0$.  Then either $c'\geq\gamma$ or $\gamma\geq c'$.  Assume that $c'\geq\gamma$, then both $m_{24}=y^bz^{\gamma}$ and $m_{45}=x^{a'}y^{\beta}z^{c'}$ strictly divide $m_{25}=x^{a'}y^bz^{c'}$.  This implies that $\sigma_{25}\not\in S_2$ by Lemma~\ref{SyzLem}.  Similarly if we assume that $\gamma\geq c'$ we have that $\sigma_{14}\not\in S_2$ which is a contradiction.  Thus if $\sigma_{14}\in S_2$ the only minimal second syzygies that can give us a $y^b$ or $z^c$ in $f_2$ are $\sigma_{12}$ and $\sigma_{13}$.

Now since the only possible minimal second syzygies that give the desired entries in $f_2$ are $\sigma_{1i}, 2\geq i\geq n$, then we have at most $n-1$ of these second syzygies.  Part (1) of the proposition says that each of these has exactly two nonzero entries from $I$ thus we can have at most $2(n-1)$ nonzero elements from $I$ as entries in $f_2$. \qed
\end{pf}

We will see later that if $I$ is generic and $n\geq5$ we achieve the maximum number of nonzero entries in $f_2$ from $I$ if and only if we achieve the maximum number of nonzero entries in $f_3$ from $I$.  This is a direct consequence of Proposition~\ref{MaxEntriesf2} and Theorem~\ref{MaxEntriesFromI}.  Refer to Example~\ref{MaxEntriesEx} to illustrate Proposition~\ref{MaxEntriesf2}.

\section{Generic Monomial Resolutions}

The focus of this section will be to discuss free resolutions of $R/I$ when $I$ is a generic monomial ideal.  Free resolutions of generic monomial ideals have been studied extensively in \cite{BPS}, \cite{EM}, \cite{MS}, and \cite{MSY}.  These resolutions have a specific structure.  In particular, $f_3$ will contain exactly three nonzero entries in each column.  Here we will discuss how we can use Buchberger graphs to represent these resolutions when $I$ is generic as shown in \cite{BPS} and \cite{MS}.  We will also discuss more specific properties of $f_3$, namely, the types of nonzero elements of $I$ we can get in $f_3$ and the maximum number of such entries.

\begin{definition}
A monomial ideal $I = \brc{m_1,\ldots,m_n}$ is \emph{generic} if whenever two distinct minimal generators $m_i$ and $m_j$ have the same positive degree in some variable, there is another minimal generator $m_k$ such that $m_k \dblbar m_{ij}$.
\end{definition}

Generic monomial ideals are defined so that $\sigma_{ij}\notin S_2$ whenever $m_i$ and $m_j$ have the same positive degree in some variable, see \cite[Theorem 6.26]{MS}.  This is also a consequence of the second syzygy lemma.  By the definition of generic, if $m_i$ and $m_j$ have the same positive degree in some variable, then there will always be a minimal generator $m_k$ such that $m_k \dblbar m_{ij}$. This implies that $m_{ik}\strdiv m_{ij}$ and $m_{jk} \strdiv m_{ij}$ by Lemma~\ref{DividesLem}.  We will now give a stronger version of Lemma~\ref{SyzLem} for generic monomial ideals.

\begin{lemma}
\label{SyzLemGen}
Let $I$ be a generic monomial ideal with minimal generating set $\{m_1, \ldots, m_n\}$, then $\sigma_{ij} \not\in S_2$ if and only if there exists a minimal generator $m_k$ such that $m_{ik}\strdiv m_{ij}$ and $m_{jk} \strdiv m_{ij}$.
\end{lemma}

\begin{pf}
 We have already proven the reverse direction in Lemma~\ref{SyzLem}.  For the forward direction we need to show that if either condition (1) or (2) in Definition~\ref{ResCons} is not satisfied then we will get the desired result.  In fact we will show that if (2) is not satisfied that this implies that (1) is also not satisfied when $I$ is generic.  First assume that (1) is not satisfied, then there exists a minimal generator $m_k$ such that $m_{ik}\strdiv m_{ij}$. We only need to show that $m_{jk}\strdiv m_{ij}$.  It can also be assumed that $m_i$ and $m_j$ do not share any positive degrees in some variable since $I$ is generic.  If $m_i$ and $m_k$ do not have the same positive degree in some variable then we are done because $m_{jk}\not=m_{ij}$.  If $m_j$ and $m_k$ do have the same positive degree in some variable, then there exists a minimal generator $m_l$ such that $m_l\dblbar m_{jk}$.  Thus $m_l\dblbar m_{ij}$.

We will now show that if (2) is not satisfied in Definition~\ref{ResCons} then (1) is not satisfied either.  If (2) is not satisfied then Lemma~\ref{SyzLem} says that there exists a minimal generator $m_k$ such that $m_k\strdiv m_{ij}$.  This implies that both $m_{ik}$ and $m_{jk}$ divide $m_{ij}$. Notice that both $m_{ik}\strdiv m_{ij}$ and $m_{jk}\strdiv m_{ij}$ would imply that (1) is not satisfied and we are done.  Without loss of generality say $m_{ik} = m_{ij}$, then at least two of these three minimal generators must have the same positive degree on the same variable. Since $I$ is generic there must exist a minimal generator $m_l$ such that $m_l$ strongly divides either $m_{ij}$, $m_{ik}$ or $m_{jk}$.  This implies that $m_l \dblbar m_{ij}$ which gives us that (1) is not satisfied.\qed
\end{pf}

Using Lemma~\ref{SyzLemGen} or \cite[Theorem 6.26]{MS} we have that all of the minimal second syzygies in $f_2$ must correspond to minimal generators $m_i$ and $m_j$ such that these generators do not have the same positive degree in some variable.  This also tells us that we will not have a choice of which second syzygies are in $S_2$.  Thus condition (2) in Definition~\ref{ResCons} is not needed to define $S_2$ for generic monomial ideals.

One of the nice properties about trivariate monomial ideals is that we can represent these ideals with three dimensional staircase diagrams.  Since we are assuming that $R/I$ artinian, this also implies that these diagrams will be bounded on all axes. In general staircase diagrams provide a template for which we can draw a graphical representation of the free resolution of $R/I$ as shown in \cite{EM} and \cite{MS}.  For a generic monomial ideal one way this can be done is by constructing the Buchberger graph for $I$.

\begin{definition}
The \emph{Buchberger graph} Buch($I$) of a monomial ideal $I = \brc{m_1, \ldots, m_n}$ has vertices  $1,\ldots, n$ and an edge $(i,j)$ whenever there is no monomial $m_k$ such that $m_k \dblbar m_{ij}$.
\end{definition}

From this definition it is clear that if $I$ is generic, then $(i,j)$ will not be an edge on Buch($I$) if $m_i$ and $m_j$ have the same positive degree in some variable.  In \cite[Lemma 6.10]{MS} it is shown that Buch($I$) is equal to the edges of the Scarf complex $\triangle_I$, which uniquely generates $Z_2$ when $I$ is generic.  Using this along with a result from \cite{BPS} we are able to give a detailed description of the free resolution of $R/I$ when $I$ is generic.

\begin{proposition}
\label{BuchProp}
Let $I = \brc{m_1, \ldots, m_n}$ be generic, then the following hold.
\begin{enumerate}
\item $(i,j) \in$ Buch($I$) if and only if $\sigma_{ij} \in S_2$.
\item Buch($I$) is a planar triangulation.
\end{enumerate}
\end{proposition}


 Part (1) follows from \cite[Lemma 6.10]{MS} and (2) is obtained from a more general result \cite[Corollary 5.5]{BPS}, which says that if $I$ is a generic monomial ideal in $r$ variables then the Scarf complex $\triangle_I$ is a regular triangulation.  Since we are only dealing with trivariate monomial ideals we can say that Buch($I$) is a planar triangulation, we refer the reader to \cite{GR} for more on planar triangulations.  In general $S_2$ is contained in the edges of Buch($I$), but Buch($I$) will not give us a minimal representation of a free resolution unless $I$ is generic.  Part (2) of Proposition~\ref{BuchProp} tells us that every column in $f_3$ contains exactly three nonzero entries.  This is not surprising since it is also known that $\triangle_I$ is a simplicial complex, \cite[Lemma 6.8]{MS}.  In Lemma~\ref{PurePower} we will give a precise description of these nonzero entries.  First, we will motivate this with an example. We will not display the first matrix, $f_1$ for some of our examples. Instead we will assume that given an ideal $I = \brc{m_1, \ldots, m_n}$ that $f_1 = \left[\begin{matrix}m_1 & \ldots & m_n\end{matrix}\right]$.

\begin{example}

Let $I = \brc{yz^2,x^5,x^3y^2,y^5,z^5,x^3z^3}$, which is generic, and the minimal free resolution for $R/I$ obtained from Definition~\ref{ResCons} is,
$$\setlength\arraycolsep{0.4mm} 0 \xrightarrow{\hspace{.05in}}  R^4 \xrightarrow{ \scriptsize \left[ \begin{array}{rrrr} y&z&0&0\\-x^2&0&y^3&0\\0&0&-x^3&0\\0&0&0&x^3\\0&-x^2&0&-z^2\\z^2&0&0&0\\0&y&0&0\\0&0&z^2&0\\0&0&0&y\end{array}\right]}  R^9 \xrightarrow{ \footnotesize \left[ \begin{array}{rrrrrrrrr}-x^5&-x^3y&-y^4&-z^3&-x^3z&0&0&0&0\\yz^2&0&0&0&0&-y^2&-z^3&0&0\\0&z^2&0&0&0&x^2&0&-y^3&0\\0&0&z^2&0&0&0&0&x^3&0\\0&0&0&y&0&0&0&0&-x^3\\0&0&0&0&y&0&x^2&0&z^2\\\end{array}\right]} R^6 \rightarrow \cdots$$

\noin The main observation here is that each column in $f_3$ contains exactly three nonzero \emph{pure power} entries.  We will show that this is true in general for minimal resolutions given by generic monomial ideals.  We first give a precise definition of the third syzygies for $R/I$ when $I$ is generic.

\end{example}

\begin{definition}
\label{ThirdSyz}
Let $I=\brc{m_1, \ldots, m_n}$ be generic, then every column in $f_3$ is given by, $$\displaystyle \mathbf{\tau}_{ijk} = \frac{m_{ijk}}{m_{ij}}e_{|\sigma_{ij}|} - \frac{m_{ijk}}{m_{ik}}e_{|\sigma_{ik}|} + \frac{m_{ijk}}{m_{jk}}e_{|\sigma_{jk}|}, \textrm{ such that } 1 \leq i < j < k \leq n,$$ and $|\sigma_{ij}|$ denotes the column number in $f_2$, in which $\sigma_{ij}$ lies.
\end{definition}
The reason that this is true is due to the fact that the following identity holds for all $1\leq i <j<k \leq n$,
$$\displaystyle \mathbf{\tau}_{ijk} = \frac{m_{ijk}}{m_{ij}}\sigma_{ij} - \frac{m_{ijk}}{m_{ik}}\sigma_{ik} + \frac{m_{ijk}}{m_{jk}}\sigma_{jk} = 0.$$

In general we will define the set of all third syzygies by $Z_3 = \sum_{i<j<k}R\mathbf{\tau}_{ijk} \subseteq R^{m+n-1}$, with minimal generating set $S_3$.  We make note that to find the minimal third syzygies we only need to find all combinations of three columns in $f_2$ so that these three columns only have nonzero entries in three distinct rows.  The $|\sigma_{ij}|$, $|\sigma_{ik}|$, and $|\sigma_{jk}|$ will be the actual column numbers of $\sigma_{ij},\sigma_{ik}$ and $\sigma_{jk}$ in $f_2$. The column numbers will also tell us which rows in $\mathbf{\tau}_{ijk}$ have nonzero entries.  It should be noted that $\mathbf{\tau}_{ijk}$ is unique up to the ordering chosen on the second syzygies.  Also since the minimal third syzygies are determined by the chosen minimal second syzygies we do not have a choice of which $\tau_{ijk}$ minimally generate $Z_3$ when $I$ is generic.

\begin{lemma}
\label{ThirdSyzLem}
Let $I$ be a monomial ideal of $R$, such that $m_i, m_j, m_k$ and $m_l$ are distinct minimal generators of $I$ with $i<j<k<l$. If $\sigma_{ij},\sigma_{ik},\sigma_{jk},\sigma_{il},\sigma_{jl},\sigma_{kl} \in S_2$ and $m_{ijk}, m_{ijl}$ and $m_{ikl}$ all strictly divide $m_{jkl}$ then $\tau_{jkl} \in \mathfrak{m}Z_3$.
\end{lemma}
\begin{pf}
First since $m_{ijk}, m_{ijl}$ and $m_{ikl}$ all strictly divide $m_{jkl}$ then,
$$\displaystyle \frac{m_{jkl}}{m_{ijk}}, \frac{m_{jkl}}{m_{ijl}}, \frac{m_{jkl}}{m_{ikl}} \in \mathfrak{m}.$$
Thus we have that,

$$\displaystyle \frac{m_{jkl}}{m_{ijk}} \tau_{ijk} - \frac{m_{jkl}}{m_{ijl}} \tau_{ijl} + \frac{m_{jkl}}{m_{ikl}} \tau_{ikl} \hspace{3in}$$ $$= \frac{m_{jkl}}{m_{ijk}}\left(\frac{m_{ijk}}{m_{ij}}e_{|\sigma_{ij}|} - \frac{m_{ijk}}{m_{ik}}e_{|\sigma_{ik}|} + \frac{m_{ijk}}{m_{jk}}e_{|\sigma_{jk}|}\right)-$$
$$\hspace{1in}\frac{m_{jkl}}{m_{ijl}}\left(\frac{m_{ijl}}{m_{ij}}e_{|\sigma_{ij}|} - \frac{m_{ijl}}{m_{il}}e_{|\sigma_{il}|} + \frac{m_{ijl}}{m_{jl}}e_{|\sigma_{jl}|}\right)+$$
$$\hspace{2in}\frac{m_{jkl}}{m_{ikl}}\left(\frac{m_{ikl}}{m_{ik}}e_{|\sigma_{ik}|} - \frac{m_{ikl}}{m_{il}}e_{|\sigma_{il}|} + \frac{m_{ikl}}{m_{kl}}e_{|\sigma_{kl}|}\right)$$

$$= \frac{m_{jkl}}{m_{jk}}e_{|\sigma_{jk}|} - \frac{m_{jkl}}{m_{jl}}e_{|\sigma_{jl}|} + \frac{m_{jkl}}{m_{kl}}e_{|\sigma_{kl}|} = \tau_{jkl}. \qed \hspace{.11in}$$
\end{pf}

In general it is not difficult to see that the set $\{\tau_{ijk} | \sigma_{ij},\sigma_{ik},\sigma_{jk} \in S_2 \textrm{ for all } i<j<k\}$ generates $Z_3$ for any monomial ideal $I$.  However, all of the second syzygies involved with an arbitrary $\tau_{ijk}$ may not be minimal.  If this occurs we must simply replace the non-minimal second syzygy with a linear combination of two minimal second syzygies that generate it, as described in Remark~\ref{LinCombTwo}. Since by definition $Z_3 \subseteq R^{m+n-1}$ the second syzygies corresponding to $\tau_{jkl}$ must be minimal to be able to say that $\tau_{jkl} \in \mathfrak{m}Z_3$. Because of this we could just assume that $\sigma_{jk}, \sigma_{jl}$ and $\sigma_{kl}$ are minimal in Lemma~\ref{ThirdSyzLem} and get the same result.

\begin{lemma}
\label{PurePower}
Let $I$ be a generic monomial ideal of $R$, then the matrix of $f_3$ in Definition~\ref{ResCons} contains only pure powers of $x$, $y$, and $z$.
\end{lemma}

The proof of this lemma is a consequence of Proposition~\ref{ThreeEntries} (which we will prove in the next section) and Definition~\ref{ThirdSyz}.

This lemma tells us that each column in $f_3$ will have exactly three nonzero entries of the form $(\pm) x^{a'}$, $(\pm) y^{b'}$ and $(\pm)z^{c'}$ where $a', b', c' > 0$. From this we can also determine exactly when we get rows in $f_3$ which contain only elements from the ideal $I$.  First we will see what conditions must be satisfied for us to get nonzero entries in $f_3$ from $I$.  By applying the previous lemma it is clear that if $x^a$, $y^b$, and $z^c$ are the minimal pure power generators of $I$,  that these are the only possible nonzero entries that we can have in $f_3$ from $I$.  So, if we wanted to have $x^a$ as an entry in $f_3$, we would need a set of three minimal generators, $\{x^a, y^{b_1}z^{c_1}, y^{b_2}z^{c_2}\}$, $b_i,c_i \geq 0$, which corresponds to a minimal third syzygy.  Notice that $\lcm{x^a, y^{b_1}z^{c_1}, y^{b_2}z^{c_2}} = x^ay^{b'}z^{c'}$ where $b' = \textrm{max}\{b_1,b_2\}$ and $c' = \textrm{max}\{c_1,c_2\}$.  Here we get $x^a$ as an entry in this third syzygy from the following computation,
$$\displaystyle \frac{\lcm{x^a, y^{b_1}z^{c_1}, y^{b_2}z^{c_2}}}{\lcm{y^{b_1}z^{c_1},y^{b_2}z^{c_2}}} = \frac{x^ay^{b'}z^{c'}}{y^{b'}z^{c'}} = x^a.$$
Now we know what is required to obtain a nonzero entry from $I$ in $f_3$.  The following lemma will show that when we get nonzero entries, from $I$ in $f_3$, then we will not get any other nonzero entries from $I$ in $f_3$, when $I$ is generic.

\begin{lemma}
\label{OnePure}
Let $I$ be a generic monomial ideal of $R$, then if the matrix of $f_3$ in Definition~\ref{ResCons} contains nonzero entries from $I$, each such entry must be the same.
\end{lemma}
\begin{pf}
To show this we will consider different possibilities for nonzero pure power generators to appear in $f_3$ and eliminate all these possibilities except for what is stated in the lemma.  First, if one column contains all three pure power generators we will see that this implies this is the only column in $f_3$, which implies that $R/I$ is Gorenstein.  For sake of contradiction suppose that there is more than one column in $f_3$. This means we must have that $I$ is minimally generated by at least four monomials.  To simplify calculation we will order our generators so that $m_1 = x^a$, $m_2=y^b$, $m_3=z^c$ and $m_4=x^{\alpha}y^{\beta}z^{\gamma}$ such that $\alpha<a$, $\beta<b$ and $\gamma<c$, where at least two of these degrees are positive.  Since $I$ is generic we know that $\tau_{123}$, $\tau_{124}$, $\tau_{134}$ and $\tau_{234}$ are all in $Z_3$ but not necessarily in $S_3$.   It is clear that $m_{124}, m_{134}$ and $m_{234}$ all strictly divide $m_{123}$, which implies that $\tau_{123} \in \mathfrak{m}Z_3$ by Lemma~\ref{ThirdSyzLem}.  Thus $\tau_{123} \not\in S_3$, which implies that if $f_3$ contains a column with all three pure power generators this must be the only column in $f_3$.

Next we will show that a single column cannot contain two different pure power generators because this would require that either more than three minimal generators had to correspond to this minimal third syzygy (which cannot happen since $I$ is generic) or $R/I$ is Gorenstein.  Without loss of generality suppose that we have $x^a$ and $y^b$ as pure power entries in the same column of $f_3$.  Then we know that we must have another generator with nonzero degrees only on $y$ and $z$ to get $x^a$ as an entry, and a generator with nonzero degrees only on $x$ and $z$ to get $y^b$ as an entry.  Since $I$ is generic we can only have one other generator corresponding to this third syzygy by Proposition~\ref{BuchProp} and Definition~\ref{ThirdSyz}.  Thus this third generator would just be $z^c$, which would imply that $R/I$ is Gorenstein by the previous argument.

We will now show that it is not possible to get two different pure power generators in two different columns of $f_3$. Suppose that $x^a$ is an entry in one column of $f_3$ and $y^b$ is an entry in another column of $f_3$. Then we must have two sets of minimal generators $\{x^a, y^{b_1}z^{c_1}, y^{b_2}z^{c_2}\}, \ b \geq b_1 > b_2 \geq 0, \ c_2 > c_1 \geq 0$ and $\{y^b, x^{a_1}z^{c_3}, x^{a_2}z^{c_4}\}, \ a \geq a_1 > a_2 \geq 0, \ c_4 > c_3 \geq 0$ that correspond respectively to these minimal third syzygies.  Here we must have that either $c_1 \geq c_3$ or $c_3 \geq c_1$, choosing either case will yield a similar contradiction.  Assume that $c_1 \geq c_3$, then $$ \lcm{x^a,x^{a_1}z^{c_3}}\strdiv \lcm{x^a, y^{b_1}z^{c_1}} \textrm{ and } \lcm{x^{a_1}z^{c_3}, y^{b_1}z^{c_1}}\strdiv \lcm{x^a, y^{b_1}z^{c_1}}.$$ Thus there is no minimal second syzygy between $x^a$ and $y^{b_1}z^{c_1}$ by Lemma~\ref{SyzLem}, which implies the minimal generators $\{x^a, y^{b_1}z^{c_1}, y^{b_2}z^{c_2}\}$ cannot correspond to a minimal third syzygy.  Therefore we cannot get two different pure power generators as entries in two different columns of $f_3$. Thus our only other option is that we may have entries from one of the pure power generators of $I$ in $f_3$. \qed
\end{pf}

 We may want to note that we could have also looked at the proof of Lemma~\ref{OnePure} by analyzing the planar graph representation of the resolution.  First to get all three pure power entries we would only be able to have one face on the graph in which the vertices would be the pure power generators, which implies $R/I$ is Gorenstein.  To get two different pure power entries in the same column, we would need at least four edges which cannot happen since the graphs associated with generic resolutions are planar triangulations.  In Section 4 we will see that this is possible when $I$ is not generic. To get two different entries in different columns we would describe two different faces such that both have a different pure power vertex and the other two vertices in the face involve only the other two variables.  But this would contradict the planarity of the graph since we would have at least two of the edges crossing at a location where there is not a vertex on the graph.
 The following example illustrates the previous two lemmas.

\begin{example}
Let $I = \brc{x^4,x^2y^2,xy^3,y^4,x^3z,z^5}$, which is generic. Then the free resolution for $R/I$ obtained from Definition~\ref{ResCons} is,
$$\setlength\arraycolsep{.6mm} 0 \xrightarrow{\hspace{.15in}} R^4 \xrightarrow{\scriptsize \left[ \begin{array}{rrrr}z&0&0&0\\-y^2&0&0&0\\0&\text{\circleme{$z^5$}}&0&0\\x&0&z^4&0\\0&-y&-x&0\\0&0&0&\text{\circleme{$z^5$}}\\0&x&0&-y\\0&0&0&x\\0&0&y^2&0\end{array}\right]} R^9 \xrightarrow{\scriptsize \left[ \begin{array}{rrrrrrrrr}-y^2&-z&0&0&0&0&0&0&0\\x^2&0&-y&-xz&-z^5&0&0&0&0\\0&0&x&0&0&-y&-z^5&0&0\\0&0&0&0&0&x&0&-z^5&0\\0&x&0&y^2&0&0&0&0&-z^4\\0&0&0&0&x^2y^2&0&xy^3&y^4&x^3\end{array}\right]} R^6 \xrightarrow{\hspace{.15in}} \cdots$$
\end{example}

In this example we have two occurrences of $z^5$ in $f_3$, which is also a minimal generator of $I$.  We obtained these two entries from the sets $\{x^2y^2,xy^3,z^5\}$ and $\{xy^3,y^4,z^5\}$ which correspond to $\tau_{2\hspace{.01in}3\hspace{.01in}6}$ and $\tau_{3\hspace{.01in}4\hspace{.01in}6}$ respectively.  Using the previous two lemmas we can now find a bound on the maximum number of nonzero entries allowed in $f_3$ which are also in $I$.

\begin{theorem}
\label{MaxEntriesFromI}
Let $I$ be a generic monomial ideal of $R$, minimally generated by $n$ elements, then the matrix of $f_3$ in Definition~\ref{ResCons} will contain at most $n-2$ nonzero entries from $I$.
\end{theorem}
\begin{pf}
Using Lemmas~\ref{PurePower} and ~\ref{OnePure} we need only show that the maximum number of entries in $f_3$ for one of the pure power generators of $I$ is $n-2$.  To do this we will show that there is no possible construction of $I$ where we can get more than $n-2$ entries in $f_3$.  We will first show exactly which combination of generators will give us entries from $I$ in $f_3$.  Without loss of generality assume $x^a$ is the entry we get in $f_3$.  Let $I$ be generic monomial ideal such that $\{x^a, y^{b_1}z^{c_1}, y^{b_2}z^{c_2}, \ldots, y^{b_r}z^{c_r}\}$ are minimal generators of $I$ with $b_r = c_1 = 0$, $b_{i-1} > b_i > 0$, and $0 < c_i < c_{i+1}$ for all $2 \leq i \leq r-1$.  Here we are not necessarily assuming that these are the only minimal generators of $I$, we are just picking out $x^a$ and all generators with nonzero degrees only on $y$ and $z$.  Notice that this still satisfies the blanket condition that $I$ is $\mathfrak{m}$-primary, since $y^{b_1}z^{c_1} = y^{b_1} = y^b$ and $y^{b_r}z^{c_r} = z^{c_r} = z^c$ by assumption.  Now we will show that the only possible third syzygies that can give us an $x^a$ in $f_3$ must involve three of these generators in the form, $\{x^a, y^{b_i}z^{c_i}, y^{b_{i+1}}z^{c_{i+1}}\}$ such that $1 \leq i \leq r-1$.
  Suppose the contrary, that the set $\{x^a,y^{b_{j}}z^{c_{j}},y^{b_i}z^{c_i}\}$ with $i,j \in \{1, \ldots, r\}$ and $|i-j| \geq 2$, corresponds to a minimal third syzygy. Then there are two cases, that $i<j$ or $i>j$.  First suppose that $i<j$, then since $b_{i} > b_{i+1} > b_j$ and $c_i < c_{i+1} < c_j$ we have that both
$$\lcm{y^{b_i}z^{c_i},y^{b_{i+1}}z^{c_{i+1}}}=y^{b_i}z^{c_{i+1}} \textrm{ and } \lcm{y^{b_{i+1}}z^{c_{i+1}},y^{b_j}z^{c_j}} = y^{b_{i+1}}z^{c_j},$$
strictly divide $\lcm{y^{b_i}z^{c_i},y^{b_j}z^{c_j}} = y^{b_i}z^{c_j}$.  Thus $\{y^{b_i}z^{c_i},y^{b_j}z^{c_j}\}$ cannot correspond to a minimal second syzygy and hence $\{x^a,y^{b_j}z^{c_j},y^{b_i}z^{c_i}\}$ cannot correspond to a minimal third syzygy.  A similar argument may be used when $i>j$.  Taking all possible sets satisfying our conditions we get,
$$\{x^a,y^{b_1}z^{c_1},y^{b_2}z^{c_2}\}, \{x^a,y^{b_2}z^{c_2},y^{b_3}z^{c_3}\}, \ldots, \{x^a,y^{b_{r-1}}z^{c_{r-1}},y^{b_{r}}z^{c_{r}}\}$$
are the only possible sets that can correspond to minimal third syzygies with $x^a$ as an entry, and there are exactly $r-1$ of these sets.  Now it is easy to see that we cannot get more than $n-2$ of these sets for an ideal generated by $n$ elements. This is because we have $r$ minimal generators with nonzero degrees only on $y$ and $z$, which gives at most $r-1$ entries with $x^a$ and we must have a minimum of $r+1$ minimal generators for $I$.  If we wanted greater than or equal to $n-1$ copies of $x^a$ in $f_3$, then $r-1 \geq n-1  \Longrightarrow r \geq n$, but this implies that $I$ must have at least $n+1$ generators, which is a contradiction.  Thus we can only have at most $n-2$ copies of $x^a$ in $f_3$, which completes our proof.\qed
\end{pf}

\begin{example}
\label{MaxEntriesEx}
Let $I = \brc{x^5,x^4y,x^2y^3,xy^4,y^5,z^5}$, which is generic. Then the free resolution for $R/I$ obtained from Definition~\ref{ResCons} is,
$$\setlength\arraycolsep{0.6mm} 0 \xrightarrow{\hspace{.15in}}  R^4 \xrightarrow{\scriptsize \left[ \begin{array}{rrrr}\text{\circleme{$z^5$}}& 0 & 0 & 0\\-y & 0 & 0 & 0\\0&\text{\circleme{$z^5$}}&0&0\\x&-y^2&0&0\\0&0&\text{\circleme{$z^5$}}&0\\0&x^2&-y&0\\0&0&0&\text{\circleme{$z^5$}}\\0&0&x&-y\\0&0&0&x\end{array}\right]} R^9 \xrightarrow{\scriptsize \left[ \begin{array}{rrrrrrrrr}-y&-z^5&0&0&0&0&0&0&0\\x&0&-y^2&-z^5&0&0&0&0&0\\0&0&x^2&0&-y&-z^5&0&0&0\\0&0&0&0&x&0&-y&-z^5&0\\0&0&0&0&0&0&x&0&-z^5\\0&x^5&0&x^4y&0&x^2y^3&0&xy^4&y^5\end{array}\right]} R^6 \xrightarrow{\hspace{.15in}} \cdots$$

Here we see that we can get $n-2$ nonzero entries from $I$ in $f_3$ and that we also have exactly $n-2$ rows in $f_3$ with nonzero entries from $I$.  Theorem~\ref{MaxEntriesFromI} shows that $n-2$ is an upper bound on the number of entries we can get from $I$ in $f_3$ when $I$ is generic.  In Example~\ref{MaxEntriesEx} we also achieve $2n-2$ nonzero entries in $f_2$ which are also elements of $I$.  In Proposition~\ref{MaxEntriesf2} we proved that this is an upper bound of such entries. If we achieve $n-2$ nonzero entries in $f_3$ from $I$, then we can see from the proof of Theorem~\ref{MaxEntriesFromI} that we have $n-1$ minimal second syzygies with nonzero entries from $I$.    Thus if $I$ is generic with $n\geq5$ we get $n-2$ nonzero entries in $f_3$ from $I$ if and only if we get $2n-2$ nonzero entries in $f_2$ from $I$.  The reason we must have that $n\geq5$ is because if $I=\brc{x^a,y^b,z^c,x^{\alpha}y^{\beta}z^{\gamma}}$ with positive degrees on all variables then we get exactly $2n-2=6$ nonzero entries from $I$ in $f_2$ but no entries from $I$ in $f_3$.  This is the only case where we can achieve the upper bound on the number of nonzero entries from $I$ and $f_2$ and not achieve the upper bound in $f_3$.
\end{example}

\section{Non-Generic Resolutions}

In this section we will consider some of the differences between the resolutions of $R/I$ when $I$ is \emph{not} generic in comparison to resolutions when $I$ is generic.  The main theorem of the section will show that when $I$ is not generic we will never have the same structure on the matrix of $f_3$ as described in the previous section for a generic ideal.  Specifically, if $I$ is not generic then there will be at least one column in $f_3$ which contains more than three nonzero entries.  We will also see the interesting nature of these results by considering examples where we relax the conditions that $I$ be a \emph{trivariate} monomial ideal that is $\mathfrak{m}-$\emph{primary}.  When we remove these restrictions on $I$ we will find that many of our results from this section do not hold.

\begin{proposition}
\label{ThreeEntries}
Let $I = \brc{m_1,\ldots,m_n}$ be a monomial ideal.  If a column in $f_3$ from Definition~\ref{ResCons} has exactly three nonzero entries then these entries must all be pure power entries.
\end{proposition}

\begin{pf}
Suppose $\tau_{ijk}$ is a minimal third syzygy.  Then by definition, $$\displaystyle \mathbf{\tau}_{ijk} = \frac{m_{ijk}}{m_{ij}}e_{|\sigma_{ij}|} - \frac{m_{ijk}}{m_{ik}}e_{|\sigma_{ik}|} + \frac{m_{ijk}}{m_{jk}}e_{|\sigma_{jk}|}, \textrm{ such that } 1 \leq i < j < k \leq n.$$  Each of the nonzero entries can be described using the following form,
$$\displaystyle \frac{m_{ijk}}{m_{ij}} = x^{\alpha_1}y^{\beta_1}z^{\gamma_1}, \ \frac{m_{ijk}}{m_{ik}} = x^{\alpha_2}y^{\beta_2}z^{\gamma_2}, \textrm{ and } \frac{m_{ijk}}{m_{jk}} = x^{\alpha_3}y^{\beta_3}z^{\gamma_3}.$$
Notice that it is not possible for $\alpha_l,\beta_l,\gamma_l > 0$ for $l\in\{1,2,3\}$, because this would imply that at least one of our generators was not minimal.  Thus we will assume that two of the powers are nonzero for one of the entries in $\tau_{ijk}$.  Without loss of generality suppose that $\displaystyle \frac{m_{ijk}}{m_{ij}} = x^{\alpha_1}y^{\beta_1}$ with $\alpha_1, \beta_1 > 0$.  This implies that $a_k > a_i, a_j$ and $b_k > b_i, b_j$.  From this we can compute the other two nonzero entries in $\tau_{ijk}$ to be
$$\displaystyle \frac{m_{ijk}}{m_{ik}} = z^{\gamma_2} \textrm{ and } \frac{m_{ijk}}{m_{jk}} = z^{\gamma_3}.$$
Now we have that $z^{\gamma_2}m_{ik} = z^{\gamma_3}m_{jk}$ which means that either $m_{ik}|m_{jk}$ or $m_{jk}|m_{ik}$.  Suppose that $m_{ik}|m_{jk}$, then $m_{ij} \strdiv m_{jk}$ since $a_k > a_i, a_j$ and $b_k > b_i, b_j$. We now have two possibilities, either, $m_{ik}=m_{jk}$ or $m_{ik}$ strictly divides $m_{jk}$.  If $m_{ik}=m_{jk}$, then $\sigma_{ik} \not\in S_2$ by Definition~\ref{ResCons} and Lemma~\ref{SyzLem}.  If $m_{ik}$ strictly divides $m_{jk}$ then $\sigma_{jk} \not\in S_2$ by Lemma~\ref{SyzLem}. Both cases contradict our assumption that $\tau_{ijk}$ is a minimal third syzygy.  Thus we can have no mixed entries in $\tau_{ijk}$ and in turn all nonzero entries in $\tau_{ijk}$ are pure powers. \qed
\end{pf}

This proposition tells us that to get nonzero mixed entries in a column of $f_3$ that we must have at least four nonzero entries in the column.  Another interesting fact that we obtain in the proof of this proposition is if we have a column in $f_3$ with exactly three nonzero entries, then we will have a pure power entry for each variable. It is not true that having four or more nonzero entries in a column of $f_3$ implies that we get a mixed entry.  For this we can just look at resolutions where $I$ is a power of the homogeneous maximal ideal of $R$.

\begin{example}
\label{NonGenEx}
Let $I = \mathfrak{m}^2$, and $J = \brc{x^3,x^2y,y^3,z^3,x^2z^2}$, neither of which are generic. Then a minimal free resolution for $R/I$ obtained from Definition~\ref{ResCons} is,
$$\setlength\arraycolsep{0.4mm} 0 \xrightarrow{\hspace{.075in}}  R^3 \xrightarrow{\scriptsize \left[ \begin{array}{rrr}z&0&0\\-y&0&0\\0&z&0\\x&-y&0\\0&x&0\\-x&0&z\\0&0&-y\\0&0&x\end{array}\right]} R^8 \xrightarrow{\scriptsize \left[ \begin{array}{rrrrrrrr}-y&-z&0&0&0&0&0&0\\x&0&-y&-z&0&0&0&0\\0&0&x&0&-z&0&0&0\\0&x&0&0&0&-y&-z&0\\0&0&0&x&y&x&0&-z\\0&0&0&0&0&0&x&y\end{array}\right]} R^6 \xrightarrow{\scriptsize\left[\begin{array}{rrrrrr}x^2&xy&y^2&xz&yz&z^2\end{array}\right]} R \xrightarrow{\hspace{.075in}} \cdots$$

\noin and a minimal resolution for $R/J$ obtained from Definition~\ref{ResCons} is,

$$\setlength\arraycolsep{0.5mm} 0 \xrightarrow{\hspace{.075in}}  R^2 \xrightarrow{\scriptsize \left[ \begin{array}{rr}z^2&0\\-y&0\\0&z^3\\x&-y^2z\\0&x^2\\0&y^3\end{array}\right]} R^6 \xrightarrow{\scriptsize\left[ \begin{array}{rrrrrr}-y&-z^2&0&0&0&0\\x&0&-y^2&-z^2&0&0\\0&0&x^2&0&-z^3&0\\0&0&0&0&y^3&-x^2\\0&x&0&y&0&z\end{array}\right]} R^5 \xrightarrow{\scriptsize \left[\begin{array}{rrrrr}x^3&x^2y&y^3&z^3&x^2z^2\end{array}\right]} R \xrightarrow{\hspace{.075in}} \cdots$$
\end{example}

As we can see in the resolution of $R/I$ all the nonzero entries in $f_3$ are pure powers but we have a column with four such entries.  In the resolution of $R/J$ we observe that there is a column in $f_3$ with four nonzero entries and that one of these entries is the mixed double, $-y^2z$.  Something else we can observe is that this entry is in row four of $f_3$ which corresponds to $\sigma_{2\hspace{.01in}6}$ in $f_2$.  This minimal second syzygy is obtained from the two minimal generators $x^2y$ and $x^2z^2$, which have the same nonzero degree on $x$.  In a generic ideal we would be guaranteed the existence of another minimal generator that would strongly divide the $[x^2y,x^2z^2]$, so that $\sigma_{2\hspace{.01in}6}$ would not be minimal.  The fact that this does not happen here is what leads to the presence of the entry $-y^2z$ in row four of $f_3$.  A consequence we can observe from this is that any column in $f_3$ that contains four or more nonzero entries must involve a minimal second syzygy, $\sigma_{ij}$ such that $m_i$ and $m_j$ have the same nonzero degree on some variable.  However the converse of this is not true.  It is possible to have a resolution so that a column in $f_3$ has exactly three nonzero entries and there is a corresponding minimal second syzygy, $\sigma_{ij}$ where $m_i$ and $m_j$ have the same nonzero degree on some variable.  This can be seen in the first column of $f_3$ for the resolution of $R/J$ in Example~\ref{NonGenEx}.  We can say that if $I$ is not generic, there will be at least two minimal generators $m_i$ and $m_j$ with the same nonzero degree on some variable such that $\sigma_{ij}\in S_2$.

\begin{lemma}
\label{NonGenSyz}
If $I = \brc{m_1,\ldots,m_n}$ is not generic then there is at least two minimal generators $m_i$ and $m_j$ with the same nonzero degree in some variable such that $\sigma_{ij} \in S_2$.
\end{lemma}

We make note that Lemma~\ref{NonGenSyz} follows from \cite[Theorem 6.26]{MS} where it lists equivalent definitions for $I$ to be generic.  Specifically, $I$ is generic if and only if the algebraic Scarf complex is a free resolution of $R/I$, and every edge $(i,j)$ of $\triangle_I$ is such that $m_i$ and $m_j$ do not have the same nonzero degree on the same variable.  We will however provide a proof of Lemma~\ref{NonGenSyz} using the construction given in Definition~\ref{ResCons} and Lemma~\ref{SyzLem}.

\begin{pf}
 To prove this we will show that either $\sigma_{ij}$ is indeed in $S_2$ or we may choose two other minimal generators satisfying our hypothesis, which will also satisfy the desired result.  Without loss of generality assume that $a_i = a_j = a'$, $b_i>b_j$ and $c_j>c_i$.  Also since $I$ is not generic we may assume that there is no minimal generator $m_k$ such that $m_k\dblbar m_{ij}$.  If $\sigma_{ij}$ does not satisfy condition (1) in Definition~\ref{ResCons} then it does not satisfy condition (2).  We will show that if $\sigma_{ij}$ does not satisfy each of the conditions in Definition~\ref{ResCons} that we will be able to find another pair of minimal generators which satisfy the desired result.

First suppose that $\sigma_{ij} \in \mathfrak{m}Z_2$.  By Lemma~\ref{SyzLem} there exists a minimal generator $m_k$ such that $m_{ik}\strdiv m_{ij}$.  This implies that either max$\{b_i,b_k\}<$ max$\{b_i,b_j\}$ or max$\{c_i,c_k\}<$ max$\{c_i,c_j\}$ since $a_i=a_j$.  But since $b_i>b_j$ by assumption we must have that max$\{c_i,c_k\}$ $<$ max$\{c_i,c_j\}$, which implies that $c_k<c_j$. We now have two possibilities either $b_i=b_k$ or $b_i>b_k$.  If $b_i=b_k$ then $m_{ij} = m_{jk}$.  We could now choose $m_i$ and $m_k$ to satisfy the conditions of our hypothesis.  If $b_i>b_k$ we may assume that $a_k = a'$, since $m_k\dblbar m_{ij}$ if $a_k<a'$.  This implies that $b_i>b_k>b_j$ and $c_j>c_k>c_i$ which in turn implies that $m_{ik}\strdiv m_{ij}$ and $m_{jk}\strdiv m_{ij}$.  We can now see that either $m_i$ and $m_k$ or $m_j$ and $m_k$ would satisfy the conditions of our hypothesis and we may change our original choice of generators.

If condition (1) of Definition~\ref{ResCons} is satisfied and (2) is not then we would have that there exists a minimal generator $m_k$ such that $m_k\strdiv m_{ij}$ which implies that $m_{ik}|m_{ij}$.  If $m_{ik}\strdiv m_{ij}$ then our argument is the same as it was above.  If $m_{ik} = m_{ij}$, then $c_k=c_j=c'$.  This implies that $a_k<a'$ and $b_i>b_k>b_j$ and hence $m_{jk} = x^{a'}y^{b_k}z^{c'}$ which strictly divides $m_{ik}$ and $m_{ij}$.  Thus we may change our choice of minimal generators to $m_j$ and $m_k$ to satisfy the conditions of our hypothesis.

In summary the previous arguments show that if $\sigma_{ij} \not\in S_2$ we can always find two other minimal generators which satisfy the conditions of our hypothesis.  Since $I$ is finitely generated we may apply all of the above arguments above inductively on the two minimal generators chosen to satisfy the conditions of our hypothesis, so that eventually we will be able to find two minimal generators $m_{i'}$ and $m_{j'}$, with the same positive degree in some variable such that there will not exist a minimal generator $m_{k'}$ such that $m_{k'}\strdiv m_{i'j'}$.  Thus $\sigma_{i'j'} \in S_2$ by Lemma~\ref{SyzLem} . \qed
\end{pf}

We can now give a description of the structure of $f_3$ from Definition~\ref{ResCons} when $I$ is not generic.

\begin{theorem}
\label{GenericConv}
If $I = \brc{m_1,\ldots,m_n}$ is not generic then there is at least one column in the matrix of $f_3$ from Definition~\ref{ResCons}, which contains more than three nonzero entries.
\end{theorem}

\begin{pf}
Since $I$ is not generic there are at least two minimal generators $m_i$ and $m_j$ with the same nonzero degree on the same variable such that $\sigma_{ij}\in S_2$ by Lemma~\ref{NonGenSyz}.  Suppose that $a_i=a_j=a'$, then $b_i > b_j$ and $c_i < c_j$ which gives that $m_{ij} = x^{a'}y^{b_i}z^{c_j}$. We notice that $\sigma_{ij}$ will correspond to two different minimal third syzygies. This is due to the fact that $\sigma_{ij}$ cannot lie on the outer boundary of the planar graph associated with the resolution of $R/I$. If this were to occur it would mean that $m_{ij}$ would only have nonzero powers on exactly two variables, which implies that $m_i$ and $m_j$ cannot have the same positive degree on some variable.  We will assume that both of the minimal third syzygies associated with $\sigma_{ij}$ contain exactly three entries and get a contradiction.  First let $\tau_{ijk}$ and $\tau_{ijl}$ be minimal third syzygies.  These two third syzygies are constructed from five minimal second syzygies, $\sigma_{ij},\sigma_{ik},\sigma_{jk},\sigma_{il}, \textrm{ and } \sigma_{jl}$.  We would only need to show that one of these second syzygies is not in $S_2$ or that either $\tau_{ijk}$ or $\tau_{ijl}$ is not in $S_3$ to get a contradiction. It is clear that $a_k,a_l \not= a'$ since this would contradict the minimality of one of our second syzygies.  We have three cases that we must consider, either $a_k > a'$ and $a_l < a'$, $a_k,a_l>a'$ or $a_k,a_l<a'$. The case of $a_k<a'$ and $a_l>a'$ would be exactly the same proof as $a_k > a'$ and $a_l < a'$.\\
(1): Let $a_k > a'$ and $a_l < a'$. Then we have the following nonzero entries for $\tau_{ijk}$,
$$\displaystyle \frac{m_{ijk}}{m_{ij}} = x^{\alpha_1}, \ \ \frac{m_{ijk}}{m_{ik}} = z^{\gamma_1}, \ \ \frac{m_{ijk}}{m_{jk}} = y^{\beta_1}.$$
Before we write the nonzero entries for $\tau_{ijl}$ we observe that since $a_l < a'$ we must have that either $b_l > b_i$ or $c_l > c_j$.  Without loss of generality say $b_l > b_i$, then the nonzero entries for $\tau_{ijl}$ will be
$$\displaystyle \frac{m_{ijl}}{m_{ij}} = y^{\beta_2}, \ \ \frac{m_{ijl}}{m_{il}} = z^{\gamma_2}, \ \ \frac{m_{ijl}}{m_{jl}} = z^{\hat{\gamma_2}}.$$
This implies that $m_{il}z^{\gamma_2} = m_{jl}z^{\hat{\gamma_2}}$ and using the same technique in the proof of Proposition~\ref{ThreeEntries} we find that either $\sigma_{il}$ or $\sigma_{jl}$ would not be minimal.  Thus $\tau_{ijl}$ is not a minimal third syzygy.\\
(2): Let $a_k,a_l>a'$.  Recall from Lemmas~\ref{DividesLem} and ~\ref{SyzLem} if there is a minimal generator $m_k\strdiv m_{ij}$ such that $m_{ik}\strdiv m_{ij}$ and $m_{jk}\strdiv m_{ij}$ then $\sigma_{ij} \not\in S_2$.  To maintain that both $\tau_{ijk}$ and $\tau_{ijl}$ are minimal we may assume that no two of these four minimal generators divides the least common multiple of the other two.  We will now construct the least common multiples for all five of the second syzygies we need here:

\hspace{1.25in} $m_{ij} = x^{a'}y^{b_i}z^{c_j}$,

\hspace{1.25in} $m_{ik} = x^{a_k}y^{\textrm{max}\{b_i,b_k\}}z^{\textrm{max}\{c_i,c_k\}}$,

\hspace{1.25in} $m_{jk} = x^{a_k}y^{\textrm{max}\{b_j,b_k\}}z^{\textrm{max}\{c_j,c_k\}}$,

\hspace{1.25in} $m_{il} = x^{a_l}y^{\textrm{max}\{b_i,b_l\}}z^{\textrm{max}\{c_i,c_l\}}$,

\hspace{1.25in} $m_{jl} = x^{a_l}y^{\textrm{max}\{b_j,b_l\}}z^{\textrm{max}\{c_j,c_l\}}$.\\
Note that it is implied that $m_k$ and $m_l$ do not divide $m_{ij}$ since $a_k,a_l>a'$.  We will now simplify the above least common multiples with the following:

\begin{enumerate}[(i)]\setlength{\itemsep}{-0pt}

\item $m_j$ does not divide $m_{ik}$ implies that $c_k < c_j$ since $b_j < \textrm{max}\{b_i,b_k\}$,
\item $m_i$ does not divide $m_{jk}$ implies that $b_k < b_i$ since $c_i < \textrm{max}\{c_j,c_k\}$,
\item $m_j$ does not divide $m_{il}$ implies that $c_l < c_j$ since $b_j < \textrm{max}\{b_i,b_l\}$,
\item $m_i$ does not divide $m_{jl}$ implies that $b_l < b_i$ since $c_i < \textrm{max}\{c_j,c_l\}$,
\item $m_l$ does not divide $m_{ik}$ implies that either $a_l > a_k$ or $c_l > \textrm{max}\{c_i,c_k\}$,
\item $m_l$ does not divide $m_{jk}$ implies that either $a_l > a_k$ or $b_l > \textrm{max}\{b_j,b_k\}$,
\item $m_k$ does not divide $m_{il}$ implies that either $a_k > a_l$ or $c_k > \textrm{max}\{c_i,c_l\}$,
\item $m_k$ does not divide $m_{jl}$ implies that either $a_k > a_l$ or $b_k > \textrm{max}\{b_j,b_l\}$.

\end{enumerate}
From this we have that:

\hspace{1.25in} $m_{ik} = x^{a_k}y^{b_i}z^{\textrm{max}\{c_i,c_k\}}$,

\hspace{1.25in} $m_{jk} = x^{a_k}y^{\textrm{max}\{b_j,b_k\}}z^{c_j}$,

\hspace{1.25in} $m_{il} = x^{a_l}y^{b_i}z^{\textrm{max}\{c_i,c_l\}}$,

\hspace{1.25in} $m_{jl} = x^{a_l}y^{\textrm{max}\{b_j,b_l\}}z^{c_j}$.\\
If $a_l > a_k$ then we have that $m_{ijk} = x^{a_k}y^{b_i}z^{c_j}, m_{ikl} = x^{a_l}y^{b_i}z^{\textrm{max}\{c_i,c_k,c_l\}}$ and $m_{jkl} = x^{a_l}y^{\textrm{max}\{b_j,b_k,b_l\}}z^{c_j}$ all strictly divide $m_{ijl} = x^{a_l}y^{b_i}z^{c_j}$, which implies that $\tau_{ijl} \in \mathfrak{m}Z_3$ by Lemma~\ref{ThirdSyzLem}.  Similarly if $a_k>a_l$ we will have that $m_{ijl}, m_{ikl}$ and $m_{jkl}$ all strictly divide $m_{ijk}$ implying that $\tau_{ijk}$ is not minimal.  We can see that $a_k \not= a_l$ because this would imply that $c_l > \textrm{max}\{c_i,c_k\}, b_l > \textrm{max}\{b_j,b_k\}, c_k > \textrm{max}\{c_i,c_l\}$ and $b_k > \textrm{max}\{b_j,b_l\}$ from conditions (v) - (viii), which cannot happen. In any of these cases we have that either $\tau_{ijk}$ or $\tau_{ikl}$ is not minimal which contradicts our assumption.\\
(3): Let $a'>a_k,a_l$.  Without loss of generality we must have that some of the exponents differ on $m_k$ and $m_l$, say $b_k > b_l$.  Then we have that $m_{jl} = x^{a'}y^{b_l}z^{c_j}$ and $m_{kl} = x^{\textrm{max}\{a_k,a_l\}}y^{b_k}z^{\textrm{max}\{c_k,c_l\}}$ both strictly divide $m_{jk} = x^{a'}y^{b_k}z^{c_j}$.  Thus by lemma~\ref{SyzLem} $\sigma_{jk}\in \mathfrak{m}Z_2$ and is not minimal which implies that $\tau_{ijk}$ is not minimal which is a contradiction.  We would see a similar result for any other choice of exponents on $m_k$ and $m_l$.

Thus we have shown for all cases that when $I$ is not generic there are at least two minimal generators $m_i$ and $m_j$ of $I$, with the same nonzero degree in some variable, so that we cannot have two minimal third syzygies $\tau_{ijk}$ and $\tau_{ijl}$.  Since these two generators must correspond with two minimal third syzygies one of these third syzygies must have more than three nonzero entries. \qed

\end{pf}

The assumption that that $I$ is a trivariate monomial ideal that is $\mathfrak{m}-$primary is crucial for this theorem.  The $\mathfrak{m}$-primary condition is what forces the minimal second syzygy $\sigma_{ij}$ to correspond to exactly two minimal third syzygies, when $m_i$ and $m_j$ have the same positive degree on the same variable.  If we remove the $\mathfrak{m}$-primary condition on $I$ then Theorem~\ref{GenericConv} does not hold.

\begin{example}
Let $I=\brc{x^4,x^3yz,x^3y^3,x^3z^3,y^3z^3}$, which is not $\mathfrak{m}$-primary.  Then the minimal free resolution for $R/I$ obtained from Definition~\ref{ResCons} is,
$$\setlength\arraycolsep{0.4mm} 0 \xrightarrow{\hspace{.15in}}  R^2 \xrightarrow{\scriptsize \left[ \begin{array}{rr}y^2&z^2\\-z&0\\0&-y\\x&0\\0&x\\0&0\end{array}\right]} R^6 \xrightarrow{\scriptsize \left[ \begin{array}{rrrrrr}-yz&-y^3&-z^3&0&0&0\\x&0&0&-y^2&-z^2&0\\0&x&0&z&0&0\\0&0&x&0&y&-y^3\\0&0&0&0&0&x^3\end{array}\right]} R^5 \xrightarrow{\hspace{.15in}} \cdots$$

\noin We observe that each column in $f_3$ contains exactly three nonzero pure power entries, but $I$ is clearly not generic.  The reason for this is that the minimal free resolution for $R/I$ here is supported on a simplicial complex.
\end{example}

This theorem gives us the converse to what we already knew about generic monomial ideals.  Specifically, that if $I$ is generic then every column in $f_3$ has exactly three nonzero pure power entries.  The following corollary gives us an alternate definition for an Artinian generic monomial ideal in three variables.

\begin{corollary}
\label{AltGenDef}
A monomial ideal $I = \brc{m_1,\ldots,m_n}$, which is primary to $\mathfrak{m}$ in $R$, is generic if and only if each column in the matrix of $f_3$ from Definition~\ref{ResCons} contains exactly three nonzero entries.
\end{corollary}

We also make an observation that these results only hold in general over \emph{trivariate} monomial ideals.  If $R=\Bbbk[x_1,\ldots,x_r]$ with $r\not=3$ and $I$ was a monomial ideal which is $\mathfrak{m}$-primary, then $R/I$ will have projective dimension $r$.  Our original question of when we get entries in the matrices of $f_i$ which are also elements of $I$ becomes significantly more complicated when $r\geq4$, even if $I$ is generic.  Although resolutions for generic monomial ideals do maintain some structure in general. Namely, the resolutions are simplicial \cite{MS}.
\begin{remark}
\label{BassRmk}

The content in this paper was originally motivated by the question of whether the first Bass number of $R/I$ is always larger than the zeroth Bass number of $R/I$.  It is known that the Betti numbers of $\omega_Q$ (the canonical module of $R/I$) are equal to the Bass numbers of $R/I$, when $R/I$ is Artinian \cite{DJL}.  An equivalent question to this is, after permissible row operations, is the number of rows in $f_3$ which contain only entries from $I$, less than or equal to $n-2$?  First we describe what we mean by permissable row operations. In general, we want to ensure that none of the rows in $f_3$ are dependant mod$I$.  We say that the $k^{th}$ row of $f_3$, denoted $r_k$, is dependent mod$I$ if, $r_k - (a_1r_1+ \cdots + a_{k-1}r_{k-1} + a_{k+1}r_{k+1} + \cdots + a_{m+n-1}r_{m+n-1}) \in IR^m$, for $a_i \in R$. When $I$ is generic the row operations do not change the number of rows in $f_3$ from Definition~\ref{ResCons}, which are contained in $I$.  This is due to the fact that each column of $f_3$ has exactly three nonzero pure power entries in each variable.  Thus if a row in $f_3$ contains an element that is not already in $I$, say an $x^{a'}$, we would be unable to get this row to be contained in $I$ by applying our row operations.  If the other two nonzero entries in this column are $y^{b'}$ and $z^{c'}$ then the only option is to take $x^{a'} - (a_1y^{b'}+a_2z^{c'})$, which cannot be in $I$ because $x^{a'}$ is not in $I$.  So for resolutions given by generic monomial ideals we only need know that the number of rows in $f_3$ from Definition~\ref{ResCons} which contain only entries from $I$ is less than or equal to $n-2$.  In Theorem ~\ref{MaxEntriesFromI} we showed that the maximum number of nonzero entries that we can get from $I$ in $f_3$ will be exactly $n-2$. Thus we can conclude that the first Bass number of $R/I$ is always larger than the zeroth Bass number of $R/I$ for resolutions given by trivariate generic monomial ideals.
\end{remark}


%

%
%
\end{document}